# DOMINO TILING IN GRID - NEW DEPENDENCE

## VALCHO MILCHEV AND TSVETELINA KARAMFILOVA

**Abstract.** This article is dedicated to domino tilings of square grids. In each of these grids domino tilings are represented using linear-recurrent sequences. For different grids are determined new dependencies.

## 1. Introduction

Domino tiling of a square grid is defined to be such a covering with $2\times1$ square tiles that each domino tile covers exactly two squares of the grid and none of the square tiles overlap.

Domino brick is defined to be parallelepiped with sizes $1\times1\times2$.

## 2. Domino Tiling of square grids – an idea of Roberto Tauraso and its development

For thoroughness we will note that the domino tilings of a square grid with sizes $2\times n$, as shown in Figure 1, are expressed by Fibonacci numbers $F_n$, where $F_0 = 1$, $F_1 = 1$ and $F_n = F_{n-1} + F_{n-2}$ when $n \geq 2$.

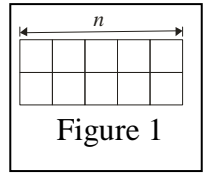

Figure 1

In 2004 Roberto Tauraso from the Tor Vergata University of Roma published an article named "A New Domino Tiling Sequence" in the specialized magazine Journal of Integer Sequences. In that article a formula for the numbers of the domino tilings of a specific grid in the shape of a right angle with sizes $2\times n\times k$ (Figure 2) is deduced using Fibonacci numbers $L_2(n,k) = F_n F_{k-1} + F_{n-1} F_k$. When $k = n-1$ Tauraso obtains the sequence with a common element $L_2(n, n-1)$, a sequence already contained in Neil Sloane's Online Encyclopedia of Integer Sequences, namely A061646. Its explicit formula is

$$L_2(n, n-1) = \frac{2}{5}\left(\frac{3+\sqrt{5}}{2}\right)^n + \frac{2}{5}\left(\frac{3-\sqrt{5}}{2}\right)^n + \frac{1}{5}(-1)^n.$$

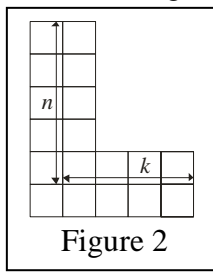

Figure 2

## Below are presented our results

We will use Tauraso's idea of a grid $L_3(2n, 2k)$ - a square grid in the form of a right angle with sizes $3 \times 2n \times 2k$, as shown in Figure 3. A grid of the defined type $L_3(2n, 2k)$ can be covered with domino tiles without them overlapping, because it is composed of two rectangular grids which by themselves can be covered with domino tiles.

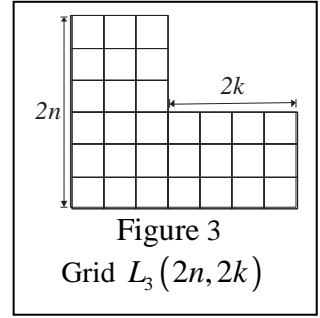
Figure 3
Grid $L_3(2n, 2k)$

**Theorem 2.1.** The number of the different ways a grid $L_3(2n, 2k)$ can be covered with domino tiles, $n \geq 2$, $k \geq 1$, is expressed by the dependence $L_3(2n, 2k) = A_n A_k + C_{n-1} B_{k-1} + B_{n-2} B_{k-1}$, where $A_n$, $B_n$ and $C_n$ express the number of domino tilings of the grids shown in Figure 4, Figure 5 and Figure 6, respectively.

**Proof.** It is known that the number of the different ways to cover the square grid $3 \times 2n$ (figure 4) with domino tiles is expressed by the formula

$$A_n = \frac{1}{6}\left[\left(3+\sqrt{3}\right)\left(2+\sqrt{3}\right)^n + \left(3-\sqrt{3}\right)\left(2-\sqrt{3}\right)^n\right] \text{ when } n \geq 1.$$

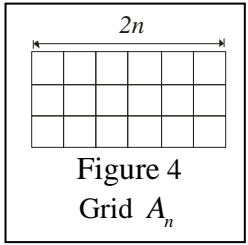
Figure 4
Grid $A_n$

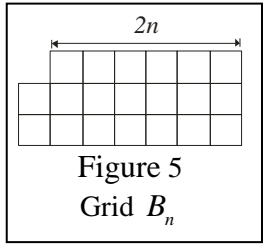
Figure 5
Grid $B_n$

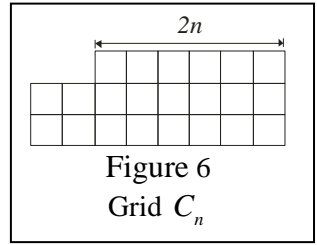
Figure 6
Grid $C_n$

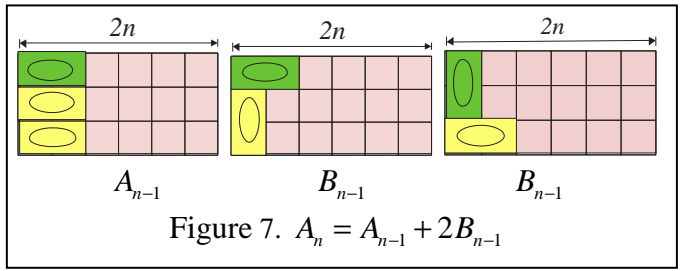
Figure 7. $A_n = A_{n-1} + 2B_{n-1}$

The formula for $A_n$ is proved using a recurrent connection of the domino tilings of grid $A_n$ and grid $B_n$. In Figure 7 the equality $A_n = A_{n-1} + 2B_{n-1}$ is shown and in Figure 8 - $B_n = A_n + B_{n-1}$.

From these two recursive equations the system of recurrent equations $A_n = 4A_{n-1} - A_{n-2}$ and $B_n = 4B_{n-1} - B_{n-2}$ is obtained. It is directly verified that $A_1 = 3$, $A_2 = 11$

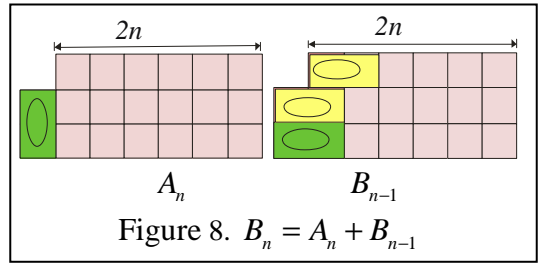
Figure 8. $B_n = A_n + B_{n-1}$

and also $B_0 = 1$, $B_1 = 4$. Explicit formulas for the two linear-recurrent sequences can be deduced.

The sequence $\{A_n\}$, where $A_n = 4A_{n-1} - A_{n-2}$ when $n \geq 3$, has characteristic equation $x^2 - 4x + 1 = 0$ with roots $x_1 = 2 + \sqrt{3}$ and $x_2 = 2 - \sqrt{3}$. We obtain formula $A_n = \frac{1}{6}\left[(3+\sqrt{3})(2+\sqrt{3})^n + (3-\sqrt{3})(2-\sqrt{3})^n\right]$ from $A_n = c_1 x_1^n + c_2 x_2^n$, $A_1 = 3$, $A_2 = 11$. For the sequence $\{B_n\}$ we have $B_0 = 1$, $B_1 = 4$, $B_n = 4B_{n-1} - B_{n-2}$, $n \geq 2$, and $B_n = \frac{1}{2\sqrt{3}}\left[(2+\sqrt{3})^{n+1} - (2-\sqrt{3})^{n+1}\right]$ when $n \geq 0$, respectively.

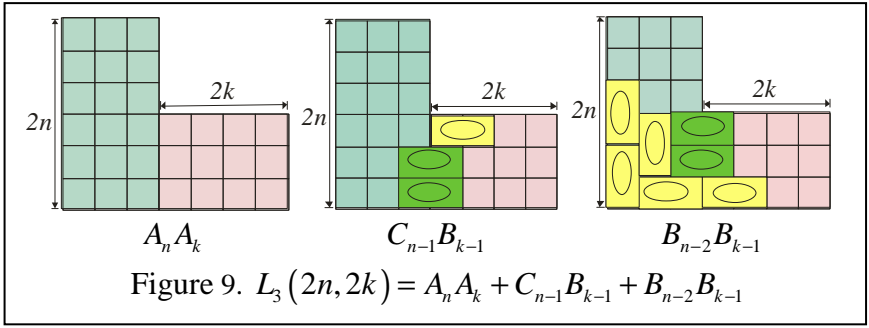

Figure 9. $L_3(2n, 2k) = A_n A_k + C_{n-1} B_{k-1} + B_{n-2} B_{k-1}$

Now look at grid $L_3(2n, 2k)$ - if we divide this grid into two rectangles – horizontal one with sizes $3 \times 2k$ and vertical one with sizes $3 \times 2n$, we conclude that it is only possible for an even number of domino tiles to cover the squares at the same time and both the parts - zero or two domino tiles. Given that the domino tiles are two, they must be adjacent. The three case with non-zero results are shown in Figure 9. Here $C_n$ is the number of the ways, the grid shown in Figure 6 can be covered with domino tiles – this is grid $3 \times (2n+2)$, from which two squares are removed, as shown in the figure. We conclude that $L_3(2n, 2k) = A_n A_k + C_{n-1} B_{k-1} + B_{n-2} B_{k-1}$.

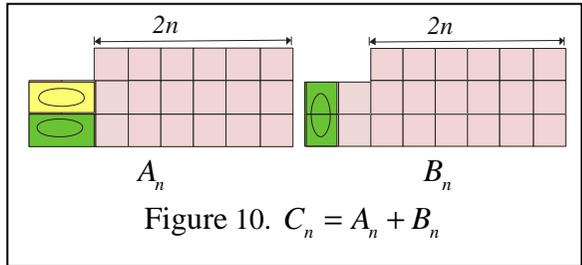

Figure 10. $C_n = A_n + B_n$

For the numbers $A_n$ and $B_n$ we have already deduced formulas. We will also deduce a formula for $C_n$. In Figure 10 the equality $C_n = A_n + B_n$ is illustrated and from it we obtain $C_0 = 2$, $C_1 = 7$, $C_n = 4C_{n-1} - C_{n-2}$ when $n \geq 2$ and also

$$C_n = \frac{1}{2}\left[(2+\sqrt{3})^{n+1} + (2-\sqrt{3})^{n+1}\right] \text{ when } n \geq 0.$$

In Table 1 the values for the sequences $\{A_n\}$, $\{B_n\}$, $\{C_n\}$ and $\{L_3(2n, 2n)\}$ are shown.

| $n$ | $A_n$ | $B_n$ | $C_n$ | $L_3(2n, 2n)$ |
|---|---|---|---|---|
| 1 | 3 | 4 | 7 | 11 |
| 2 | 11 | 15 | 26 | 153 |
| 3 | 41 | 56 | 97 | 2 131 |
| 4 | 153 | 209 | 362 | 29 681 |
| 5 | 571 | 780 | 1 351 | 413 403 |
| 6 | 2 131 | 2 911 | 5 042 | 5 757 961 |
| 7 | 7 953 | 10 864 | 18 817 | 80 198 051 |
| 8 | 29 681 | 40 545 | 70 226 | 1 117 014 753 |
| 9 | 110 771 | 151 316 | 262 087 | 15 558 008 491 |
| 10 | 413 403 | 564 719 | 978 122 | 216 695 104 121 |

Table 1

**Theorem 2.1. (Crux Mathematicorum, Vol. 42(3), 2016, Problem 4128, Proposed by Valcho Milchev and Tsvetelina Karamfilova)** Let $A_n$ be the number of domino tilings of grid $3 \times 2n$. Let $L_3(2n, 2n)$ be the number of domino tilings of grid from the type, shown in Figure 3. It is composed of two grids: $3 \times 2n$ and $3 \times 2k$, where $n \geq 2$ and $k \geq 1$ are natural numbers. Prove that $L_3(2n, 2n) = A_{2n}$.

**Proof.** When $n = k$ we have the recurrent equality
$$L_3(2n, 2n) = A_n A_n + C_{n-1} B_{n-1} + B_{n-2} B_{n-1}.$$
Now we obtain
$$A_n A_n = \frac{1}{6}\left[\left(2+\sqrt{3}\right)^{2n+1} + 2 + \left(2-\sqrt{3}\right)^{2n+1}\right];$$
$$C_{n-1} B_{n-1} = \frac{\sqrt{3}}{12}\left[\left(2+\sqrt{3}\right)^{2n} - \left(2-\sqrt{3}\right)^{2n}\right];$$
$$B_{n-2} B_{n-1} = \frac{1}{12}\left[\left(2-\sqrt{3}\right)\left(2+\sqrt{3}\right)^{2n} - 4 + \left(2+\sqrt{3}\right)\left(2-\sqrt{3}\right)^{2n}\right]$$
and respectively
$$L_3(2n, 2n) = \frac{1}{6}\left[\left(3+\sqrt{3}\right)\left(7+4\sqrt{3}\right)^n + \left(3-\sqrt{3}\right)\left(7-4\sqrt{3}\right)^n\right].$$
Therefore
$$A_{2n} = \frac{1}{6}\left[\left(3+\sqrt{3}\right)\left(2+\sqrt{3}\right)^{2n} + \left(3-\sqrt{3}\right)\left(2-\sqrt{3}\right)^{2n}\right] =$$
$$\frac{1}{6}\left[\left(3+\sqrt{3}\right)\left(7+4\sqrt{3}\right)^n + \left(3-\sqrt{3}\right)\left(7-4\sqrt{3}\right)^n\right] = L_3(2n, 2n).$$

From the resulting formula it is concluded that the sequence $L_3(2n, 2n)$ is linear recurrent - it has the characteristic equation $x^2 - 14x + 1 = 0$ and the linear recurrent dependence $L_3(n) = 14L_3(n-1) - L_3(n-2)$, $n \geq 3$, briefly noted $L_3(n) = L_3(2n, 2n)$.

**Note.** The sequences $\{A_n\}$, $\{B_n\}$, $\{C_n\}$ and $\{L_3(2n, 2n)\}$ are contained in The Online Encyclopedia of Integer Sequences with their respective numbers A001835, A001353, A001075 and A122769. For the sequence A061646 it is stated that it expresses the domino tilings of the grid $L_2(n, n-1)$ from the result of Roberto Tauraso, but for the sequence A122769 it is not stated that it expresses the domino tilings of the grid $\{L_3(2n, 2n)\}$.

### 3. Building a tower with sizes $2 \times 2 \times n$ with domino bricks

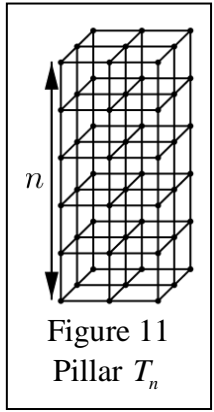

Figure 11
Pillar $T_n$

**Theorem 3.1.** Let $T_n$ be the number of the different ways a tower in the shape of a parallelepiped with sizes $2 \times 2 \times n$ can be covered with domino bricks with sizes $1 \times 1 \times 2$ (Figure 11), where $n$ is whole positive number. Then

$$T_n = \frac{1}{6}(2+\sqrt{3})^{n+1} + \frac{1}{6}(2-\sqrt{3})^{n+1} + \frac{1}{3}(-1)^n \text{ when } n \geq 1.$$

**Proof.** It is not difficult to count $T_1 = 2$, $T_2 = 9$. Let's start building our tower. In Figure 12 are shown all the possible ways of placing the first domino bricks – two horizontal bricks, four vertical bricks and a combination of two horizontal and two vertical bricks. Hence we obtain the equality $T_n = 2T_{n-1} + T_{n-2} + 4M_{n-1}$, $n \geq 3$. Here $M_n$ is the number of the ways of covering the grid, shown in Figure 13, with domino bricks.

Now we look at tower $M_n$. One can easily check that $M_1 = 1$, $M_2 = 3$. In Figure 14 are shown the possible ways of placing the first domino bricks – one horizontal or two vertical. Now we have a second recurrent relation $M_n = T_{n-1} + M_{n-1}$.

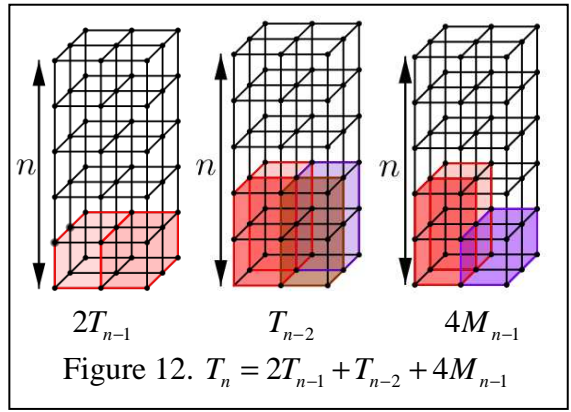

$2T_{n-1}$     $T_{n-2}$     $4M_{n-1}$
Figure 12. $T_n = 2T_{n-1} + T_{n-2} + 4M_{n-1}$

We obtain two recurrent equations:
$T_n = 2T_{n-1} + T_{n-2} + 4M_{n-1}$, $M_n = T_{n-1} + M_{n-1}$.
From this system we find
$T_n = 3T_{n-1} + 3T_{n-2} - T_{n-3}$, $M_n = 3M_{n-1} + 3M_{n-2} - M_{n-3}$, $n \geq 3$.

From the deduced dependencies it follows that the two linear-recursive sequences have a characteristic equation $x^3 - 3x^2 - 3x + 1 = 0$. Its roots are $x_1 = 2 + \sqrt{3}$, $x_2 = 2 - \sqrt{3}$ and $x_3 = -1$. Using $T_1 = 2$, $T_2 = 9$ and $T_3 = 32$ we obtain the formula

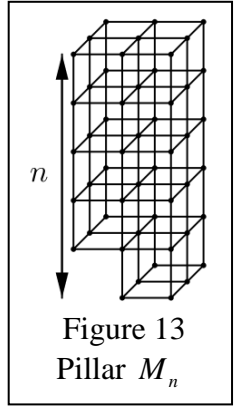

Figure 13 Pillar $M_n$

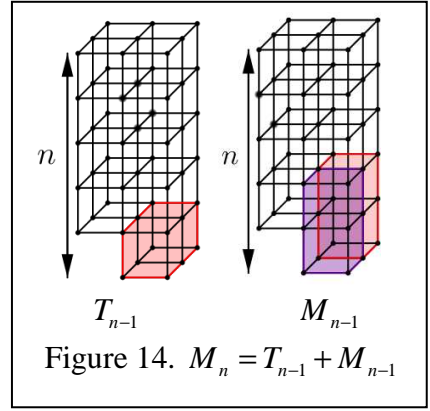

Figure 14. $M_n = T_{n-1} + M_{n-1}$

$$T_n = \frac{1}{6}\left(2+\sqrt{3}\right)^{n+1} + \frac{1}{6}\left(2-\sqrt{3}\right)^{n+1} + \frac{1}{3}(-1)^n.$$

**Theorem 3.2.** Let $T_n$ be the number of the different ways a tower in the shape of parallelepiped with sizes $2 \times 2 \times n$ can be covered with domino bricks with sizes $1 \times 1 \times 2$ (Figure 11), where $n$ is a positive whole number. Let $A_n$ be the number of the domino tilings of the grid, shown in Figure 4 and $B_n$ be the number of the domino tilings of the grid, shown in Figure 5. Then $T_{2n} = A_n^2$, $T_{2n+1} = 2B_n^2$.

| $n$ | $A_n$ | $B_n$ | $T_n$ |
|---|---|---|---|
| 1 | 3 | 4 | 2 |
| 2 | 11 | 15 | 9 |
| 3 | 41 | 56 | 32 |
| 4 | 153 | 209 | 121 |
| 5 | 571 | 780 | 450 |
| 6 | 2 131 | 2 911 | 1681 |
| 7 | 7 953 | 10 864 | 6272 |
| 8 | 29 681 | 40 545 | 23409 |
| 9 | 110 771 | 151 316 | 87362 |
| 10 | 413 403 | 564 719 | 326041 |

Table 2

**Proof.** From Table 2 it can be noted that the obtained values of $T_n$ are satisfied by $T_{2n} = A_n^2$, $T_{2n+1} = 2B_n^2$. In the general case we have

$$A_n^2 = \left\{\frac{1}{6}\left[\left(3+\sqrt{3}\right)\left(2+\sqrt{3}\right)^n + \left(3-\sqrt{3}\right)\left(2-\sqrt{3}\right)^n\right]\right\}^2 =$$

$$= \frac{1}{6}\left(2+\sqrt{3}\right)^{2n+1} + \frac{1}{6}\left(2-\sqrt{3}\right)^{2n+1} + \frac{1}{3} = T_{2n},$$

$$2B_n^2 = 2\left\{\frac{1}{2\sqrt{3}}\left[\left(2+\sqrt{3}\right)^{n+1} - \left(2-\sqrt{3}\right)^{n+1}\right]\right\}^2 =$$

$$= \frac{1}{6}\left(2+\sqrt{3}\right)^{2n+2} + \frac{1}{6}\left(2-\sqrt{3}\right)^{2n+2} - \frac{1}{3} = T_{2n+1}.$$

**Note**. The sequence $\{T_n\}$ is contained in The Online Encyclopedia of Integer Sequences with number A006253, but the connection of the domino tilings of a parallelepiped with sizes $2\times 2\times n$ and the sequences $\{A_n\}$ and $\{B_n\}$ is not pointed out.

Valcho Milchev
"Petko Rachov Slaveikov" Secondary School
Kardzhali, Bulgaria

Tsvetelina Karamfilova
"Petko Rachov Slaveikov" Secondary School
Kardzhali, Bulgaria

Valcho Milchev, teacher, Petko Rachov Slaveikov Secondary School
e-mail:
milchev.vi@gmail.com
milchev_v@abv.bg